\documentclass[12pt, reqno]{amsart}
\usepackage{graphicx} 
\usepackage{fancyhdr}

\title{Outline}
\author{Seljon Akhmedli}

\title[Hidden frequencies in topological dynamical systems]{Hidden frequencies in topological dynamical systems}

\thanks{The author was partially supported by NSF grant DMS-2136217.}

\usepackage{amsmath, amssymb, mathrsfs}
\usepackage{mathpazo}
\usepackage{latexsym,array,delarray,amsthm,amssymb,epsfig,setspace,tikz,amsmath,enumerate,mathrsfs,graphicx, mathtools,cite}

\usepackage{tabularx}
\usepackage{tikz-cd}
\usepackage{tikz}
\usepackage{floatrow}
\usepackage{hyperref}
\usetikzlibrary{positioning}
\usepackage{comment}
\usepackage{amsthm}
\usepackage{arydshln}

\usepackage{geometry}

\theoremstyle{definition}

\newcommand{\supp}{\mathop{\rm supp}}

\makeindex

\def\lnfrac#1#2{\raise.7ex \hbox{\Small $#1$}
  \kern-.15em/\kern-.15em  \lower.2ex \hbox{\Small $#2$}}

\theoremstyle{plain}
\newtheorem{theorem}{Theorem}[section]

\newtheorem{lemma}[theorem]{Lemma}
\newtheorem{proposition}[theorem]{Proposition}

\theoremstyle{definition}

\newtheorem{remark}[theorem]{Remark}

\def\thechsec{{\thechapter}.{\thesection}}
\numberwithin{equation}{section}

\usepackage{setspace}

\begin{document}
\pagenumbering{arabic}

\def\thechsec{{\thechapter}.{\thesection}}
\numberwithin{equation}{section}
\begin{abstract}
 We show for certain skew products built over any infinite minimal nilsystem, the frequencies present in the nilsystem are not always enough to govern recurrence for the whole system. This strengthens and generalizes a theorem of Glasscock, Koutsogiannis, and Richter.
\end{abstract}

\maketitle

\section{Introduction}

One of the fundamental themes arising in dynamical systems is recurrence, studying when points and sets in a system return arbitrarily close to themselves. In this article, we focus on topological recurrence for minimal dynamical systems $(X,T)$ where $X$ is a compact metric space and $T$ is a minimal homeomorphism of $X$. A system is \textit{minimal} if the orbit of every point under $T$ is dense. We often find it more useful to endow $X$ with the Borel sigma-algebra $\mathcal{B}$ equipped with a $T$-invariant Borel probability measure and look at the measure-preserving system $(X,\mathcal{B},\mu,T)$. We omit the sigma-algebra throughout and write $(X,\mu,T)$.

A \textit{set of recurrence $R \subseteq \mathbb{N}$ for a system $(X,T)$} is such that some point $x\in X$ returns arbitrarily close to itself at times belonging to $R$, namely, $\displaystyle \inf_{n\in R}d(x,T^nx)=0$. For example, $\mathbb{N}$ is easily seen to be a set of recurrence for any system by compactness. 
More interesting examples include the set of shifted primes $\mathbb{P}+1$ and the set of squares $\{n^2\colon n\in \mathbb{N}\}$ (see \cite{ Furstenberg squares, Sarkozy squares, Sarkozy}). These are sets of recurrence for all systems. On the other hand, given $1,\alpha, \beta$ rationally independent over $\mathbb{Q}$, we can quickly construct a set of recurrence for rotation by $\alpha$ on $\mathbb{T}$ which is not a set of recurrence for rotation by $\beta$ on $\mathbb{T}$ using only minimality of rotation by $(\alpha,\beta)$ in $\mathbb{T}^2$. 

It is a long-standing open question, first appearing in a paper of Veech \cite{Veech} and later popularized by Katznelson \cite{Katznelson}, whether every set of recurrence for all finite dimensional toral rotations is a set of recurrence for all systems. The question was answered in the positive for nilsystems in \cite{Host Kra Maass} and certain toral skew products in \cite{Florian et al}. An interesting parallel question one is tempted to ask is: given a minimal system $(X,T)$, \textit{which} toral rotations suffice for governing recurrence in $X$? Perhaps a natural candidate would be the \textit{maximal equicontinuous factor} of $X$ (denoted by $X_1$). Indeed, given a nilsystem, a set of recurrence for its $X_1$ factor always lifts to being a set of recurrence for the system itself following the proof in \cite[Theorem 4.1]{Host Kra Maass}. However, in \cite{Florian et al}, it is shown this behavior already fails for a skew product on $\mathbb{T}^{2}$ of the form $$T(x,y)=(x+\alpha,y+f(x))$$ for  \textit{certain} choices of $\alpha \notin \mathbb{Q}$ which are of Liouville type and $f\colon \mathbb{T} \rightarrow \mathbb{T}$ continuous where the $X_1$ factor of $(\mathbb{T}^2,T)$ is $(\mathbb{T},T_{\alpha})$ where $T_\alpha\colon \mathbb{T} \rightarrow \mathbb{T}$ is given by $x \mapsto x+\alpha$. In other words, the extension by $f$ produces hidden frequencies not captured by the maximal equicontinuous factor. We show the same phenomenon occurs in fact \textit{for all} $\alpha \notin \mathbb{Q}$ and more generally for any infinite minimal nilsystem. To obtain the right generalization, we draw on the higher order analogs of the $X_1$ factor, the maximal $k$-step nilfactors $X_k$, constructed in \cite{topological nilfactors}.

\begin{theorem}\label{main theorem}
    Let $(X,T)$ be an infinite minimal nilsystem. There exists a continuous map $f\colon X \rightarrow \mathbb{T}$ for which the system $(X\times \mathbb{T},T_f)$
    $$T_f(x,y)=(Tx,y+f(x))$$
    satisfies the following:
    \begin{enumerate}
        \item $(X\times \mathbb{T},T_f)$ is minimal and $X_k(T)\overset{t}\cong X_k(T_f)$ for all $k \in \mathbb{N}$
        \item there exists a set of recurrence for $(X,T)$ which is not a set of recurrence for $(X\times \mathbb{T}, T_f).$
    \end{enumerate}
\end{theorem}

\begin{remark}
    We later prove a more detailed statement in Section \ref{proof of main theorem} by considering the Haar measure on $X$ and the product measure $\mu \times m_{\mathbb{T}}$ on $X \times \mathbb{T}$ where $m_{\mathbb{T}}$ is the Haar measure of $\mathbb{T}$ to obtain measure-preserving systems $(X,\mu,T)$ and $(X \times \mathbb{T},\mu \times m_{\mathbb{T}},T_f)$. As an addition to property 1, we prove that $(X \times \mathbb{T},\mu\times m_{\mathbb{T}},T_f)$ is strictly ergodic and $Z_1(T_f)$ is not measurably isomorphic to $X_1(T_f)$ (see Section \ref{definitions} and Section \ref{nilsystems} for definitions and notations).
\end{remark}

Throughout this article, let $\mathbb{T}=\mathbb{R}/\mathbb{Z}$ be the one-dimensional torus. To construct the desired $f$ in Theorem \ref{main theorem}, we draw inspiration from an example of Furstenberg \cite{Furstenberg} of a minimal but not uniquely ergodic system on $\mathbb{T}^{2}$ whose dynamics is given by $(x,y)\mapsto(x+\alpha,y+f(x))$ where $\alpha \notin \mathbb{Q}$ and $f$ is a measurable but not continuous coboundary over $(\mathbb{T},m_{\mathbb{T}},T_{\alpha})$. 
In general, for a system $(X,\mu,T)$, we say $h\colon X \rightarrow \mathbb{T}$ is a \textit{measurable but not continuous coboundary over $(X,\mu,T)$} if $h$ is continuous and there exists a measurable solution $ H\colon X \rightarrow \mathbb{T}$ to  
    \begin{equation} \label{eq:coboundary}
        h=H \circ T - H
    \end{equation} but there exists no continuous solution to (\ref{eq:coboundary}). We say $h\colon X \rightarrow \mathbb{T}$ is a \textit{continuous coboundary} if there exists a continuous solution $H$ to (\ref{eq:coboundary}). The same definitions apply for $h\colon X \rightarrow \mathbb{R}$. Note all solutions to equation (\ref{eq:coboundary}) are unique up to an additive constant if $(X,\mu,T)$ is ergodic. Indeed, if $h =f\circ T-f=g\circ T-g$ for some $f,g\colon X \rightarrow \mathbb{T}$, then $f-g$ is $T$-invariant, hence constant almost everywhere by ergodicity. Section \ref{coboundaries over group rotations} is devoted to building measurable but not continuous coboundaries over any minimal rotation on an infinite compact abelian group, and more generally over any minimal system with an infinite $X_1$ factor. 

    \medskip

     \noindent\textbf{Acknowledgments.}
The author thanks Bryna Kra for her continued guidance throughout the project and many helpful suggestions on different drafts. The author would also like to thank Vitaly Bergelson, Amadeus Maldonado, and Florian Richter for fruitful discussions. 

\noindent \textbf{Use of technology:} No AI tools were used to produce any of the ideas and results. The author asked GPT 5.6 Luna some auxiliary questions related to Lemma \ref{getting a torus valued coboundary} which was partly helpful.

\section{Preliminaries}

\subsection{Definitions and notations}\label{definitions}
Throughout, $\mathbb{N}$ denotes the set $\{1,2, 3, \dots \}$. 
For a system $(X,\mu,T)$, we say $(Y,\nu,S)$ is a \textit{measurable factor} of $(X,\mu,T)$ if there exists a measurable map $\pi\colon X \rightarrow Y$, \textit{the measurable factor map}, such that $(S \circ \pi)(x)= (\pi \circ T)(x)$ for $\mu$-almost every $x \in X$ and $\pi_{*} \mu=\nu.$ Furthermore, if $\pi$ is a bijection with $\pi^{-1}$ measurable, then it is a \textit{measurable isomorphism}. We say $(Y, \nu, S)$ is a \textit{topological factor} of $(X, \mu, T)$ if there exists a continuous and surjective map $\rho\colon X \rightarrow Y$, \textit{the topological factor map}, such that $(S \circ \rho)(x)= (\rho \circ T)(x)$ for all $x\in X$ and $\rho_{*}  \mu=\nu.$ If $\rho\colon X \rightarrow Y$ is a homeomorphism, then it is a \textit{topological isomorphism}. We use the notation $\overset{m}\cong$ and $\overset{t}\cong$ to denote measurable and topological isomorphisms, respectively. Given a system $(X,\mu,T)$ and $f\in C(X)$, we at times use the notation $f_n$ to denote $\displaystyle\sum_{i=0}^{n-1}f\circ T^i$.

\subsection{Nilsystems and nilfactors}\label{nilsystems}
Given a $k$-step nilpotent Lie group $G$ and $\Gamma \leq G$ a discrete cocompact subgroup, $X=G/\Gamma$ is a \textit{$k$-step nilmanifold}. Now fix $g \in G$. Then $T_g\colon X \rightarrow X$ given by $x \mapsto  g \cdot x$ is a \textit{nilrotation.} The Haar measure on $X$ is the unique measure invariant under all nilrotations and we denote it by $m_X$. A system $(X,m_X,T)$ where $X$ is a $k$-step nilmanifold and $T$ is a nilrotation is a \textit{$k$-step nilsystem}. (All compact abelian groups are assumed to be metrizable and monothetic throughout.) A measurable inverse limit of $k$-step nilsystems is a \textit{measurable $k$-step pronilsystem} and a topological inverse limit of $k$-step nilsystems is a \textit{topological $k$-step pronilsystem}.

A system $(X,\mu,T)$ is \textit{minimal} if it has no nonempty proper closed $T$-invariant subsets of $X$ and it is \textit{ergodic} if every measurable $T$-invariant subset of $X$ has either zero or full measure. We say $(Y,\nu,S)$ is a \textit{maximal factor of $(X,\mu,T)$ with property $\mathcal{P}$} if every factor of $X$ with Property $\mathcal{P}$ is a factor of $Y$. For an ergodic system $(X, \mu, T)$, \textit{the measurable $k$-step nilfactor of $X$}, denoted by $(Z_{k}(X),m_{Z_k(X)},T)$, is the maximal measurable factor that is a measurable $k$-step pronilsystem, \textit{maximal} being with respect to measurable factor maps \cite{Host Kra}.
Similarly, if we have a minimal system $(X,T)$, then \textit{the topological $k$-step nilfactor of $X$}, denoted by $(X_{k}(X),T)$, is the maximal topological factor that is a topological $k$-step pronilsystem, where \textit{maximal} is now with respect to topological factors \cite{topological nilfactors}. For notational simplicity, we at times use $Z_k(T)$ and $X_k(T)$ to denote the measurable and topological $k$-step nilfactors, and if the transformation is clear from context, we write $Z_k$ and $X_k$. Since $Z_k$ is the \textit{maximal} measurable $k$-step nilfactor, there exists a measurable factor map from $Z_k \rightarrow X_k$, $k\in \mathbb{N}$.

The factor $X_k$ is naturally endowed with Haar measure because for pronilsystems, minimality is equivalent to unique ergodicity (only one invariant measure exists), the unique measure being Haar (see \cite{topological nilfactors}). In this way, we can view $X_k$ not only as a topological factor of our system but also as a measurable factor $(X_k(X),m_{X_k(X)},T)$. Systems that are both minimal and uniquely ergodic are \textit{strictly ergodic}.

\subsection{Eigenfunctions}
    Given a minimal and ergodic system $(X,\mu,T)$, we say $f\colon X \rightarrow \mathbb{T}$ is a \textit{measurable eigenfunction} if $f \in L^2(X)$ and there exists $c \in \mathbb{T}$ such that $f \circ T=f+c$. If $f$ is continuous, then it is a \textit{continuous eigenfunction}. An equivalent description of the $Z_1$ and $X_1$ factors which is often more useful in practice is that $L^2(Z_1)$ is the closed linear span in $\mathbb{C}$ of the measurable eigenfunctions of $X$ while $C(X_1)$ is the uniformly closed linear span in $\mathbb{C}$ of the continuous eigenfunctions of $X$. We write $\mathrm{Eig}(X,T)$ to denote the set of eigenvalues corresponding to measurable eigenfunctions of $(X,\mu,T)$. Also, we use the notation $e(t)=e^{2\pi i t}$ for all $t \in \mathbb{R}$.

\section{Furstenberg coboundaries over compact abelian group rotations}\label{coboundaries over group rotations}

\begin{lemma}\label{coboundary over group rotation}
    Let $(G, m_G,T_\alpha)$ be a minimal rotation on an infinite compact abelian group such that $\widehat{G}$ is a torsion group. Then there exists $h \colon G \rightarrow \mathbb{R}$ such that $rh$ is a measurable but not continuous coboundary over $(G,m_G,T_\alpha)$ for all nonzero $r\in \mathbb{R}$.
\end{lemma}

\begin{proof}
    Let $G=\overline{\langle\alpha\rangle}$. 
    Then there exists an infinite sequence $\gamma_n \in \widehat{G}$ such that $\gamma_n(\alpha) \rightarrow 1$. (Consider the subgroup $S=\{\phi(\alpha)\colon\phi  \in \widehat{G}\}$. Since $G=\overline{\langle\alpha\rangle}$, the map $\phi\mapsto\phi(\alpha)$, $\phi\in \widehat{G}$, is injective. It follows $S$ is infinite. By compactness, there exists a sequence $\phi_n$ such that $\phi_n(\alpha) \rightarrow \beta$ with $\phi_n(\alpha)$ distinct; thus $\phi_{n+1}\phi_n^{-1}(\alpha)\rightarrow 1$. If the sequence $\{\phi_{n+1}\phi_{n}^{-1}\}$ were  finite, then there would exist $N\in \mathbb{N}$ such that for all $n\geq N$, $\phi_{n+1}(\alpha)=\phi_n(\alpha)$ giving a contradiction. Hence $\{\phi_{n+1}\phi_n^{-1}\}$ is infinite and distinct up to passing to a subsequence.) Choose a subsequence $\gamma_{n_k}$ such that $|\gamma_{n_k}(\alpha)-1|<1/k$.
    Define $h\colon G \rightarrow \mathbb{R}$ by 
    $$h(g)=\sum_{k\neq 0}\frac{\gamma_{n_k}(\alpha)-1}{|k|}\gamma_{n_k}(g)$$
    where $\gamma_{n_{-k}}:=\overline{\gamma_{n_k}}$. 
    We write $h(g)=H(g+\alpha)-H(g)$ where $$H(g)=\sum_{k\neq 0}\frac{1}{|k|}\gamma_{n_k}(g)\;\;\mathrm{in\;}  L^2(G).$$ Notice $h$ is continuous because $\displaystyle \sum_{k\neq 0}|\widehat{h}(\gamma_{n_k})|<\infty$ and $H\in L^2(G)$ because $\displaystyle \sum_{k \neq 0} |\widehat{H}(\gamma_{n_k})|^2<\infty$. Let $\{F_k\}$ be a Folner sequence for $\widehat{G}$ such that $F_k \subseteq F_{k+1}$ and $F_k\leq \widehat{G}$ for each $k \in \mathbb{N}$. 
    Define $$\Phi_{F_k}(g)=\sum_{\gamma \in F_k}\gamma(g).$$ Then $\int \Phi_{F_k}=\widehat{\Phi_{F_k}}(\textbf{1})=1$. Notice $$C_k:=\{g \in G\colon \gamma(g)=1 \mathrm{\;for \; every \;} \gamma\in F_k\}=\supp(\Phi_{F_k}).$$ (If $g \in C_k$, then $\Phi_{F_k}(g)=|F_k|$. If $g \notin C_k$, there exists $\gamma' \in F_k$ with $\gamma'(g)\neq 1$. Since $F_k$ is a subgroup, $\gamma'F_k=F_k$ and we have $$\sum_{\gamma \in F_k}  \gamma(g)=\sum_{\gamma \in F_k} (\gamma' \gamma)(g)=\gamma'(g)\sum_{\gamma \in F_k}\gamma(g).$$ Then $\displaystyle(1-\gamma'(g))\sum_{\gamma \in F_k}\gamma(g)=0$, and since $\gamma'(g)\neq 1$, we obtain $\displaystyle\sum_{\gamma \in F_k}\gamma(g)=0$.) In particular, $\Phi_{F_k}(g)\geq 0$. 
    Since $\displaystyle\widehat{G}=\bigcup_{k}F_k$ and characters separate points of $G$, we see $\displaystyle \bigcap_{k}C_k = e$. Also, $C_{k+1}\subseteq C_k$ so we have a sequence of nested closed sets shrinking to $e$.
    Then for each $\Phi_{F_k}$ there exists an open set $U_k$ with $C_k \subseteq U_k$ and $U_k \rightarrow e$. Thus, $\Phi_{F_k}$ satisfies properties (i) and (ii) of \cite[Proposition 2.44]{Folland}. Furthermore, $$H * \Phi_{F_k}(g)=\sum_{\gamma \in \widehat{G}}\widehat{H}(\gamma) \widehat{\Phi_{F_k}}(\gamma) \gamma(g)$$ and using that $\widehat{\Phi_{F_k}}(\gamma)=1$ whenever $\gamma \in F_k$ and is 0 otherwise, we see $H*\Phi_{F_k}$ diverges at $e$. By \cite[Proposition 2.44]{Folland}, $H$ cannot be continuous. Any solution to the coboundary equation for $h$ is of the form $H+c$ where $c$ is a constant and $\widehat{(H+c)}(\gamma)=\widehat{H}(\gamma)$ for all $\gamma$ nontrivial. Thus, $H+c$ is also not continuous. Finally, for all $r\in \mathbb{R}\setminus\{0\}$, $rh$, is a measurable but not continuous coboundary because $rH*\phi_{F_k}$ diverges at $e$. \end{proof}


\begin{lemma}\label{coboundary lifts}
    Let $(X,\mu ,T)$ be minimal, $\pi\colon X \rightarrow Y$ be a topological factor map, $f\colon Y \rightarrow \mathbb{R}$ continuous, and $f\circ \pi \colon X \rightarrow \mathbb{R}$ be a continuous coboundary over $(X,\mu,T)$. Then $f$ is a continuous coboundary over $(Y,\nu,S)$.
\end{lemma}

\begin{proof}
    Applying the Gottschalk-Hedlund theorem \cite[Theorem 4.11]{Gottschalk Hedlund}, we equivalently show the ergodic sums $\displaystyle \sum_{i=1}^{n} f(S^iy)$ are uniformly bounded. Since $\pi$ is a topological factor map, for every $y \in Y$ there exists $x \in X$ such that $$\sum_{i=1}^{n} f(S^iy)=\sum_{i=1}^{n}f(S^{i}\pi(x))=\sum_{i=1}^{n}(f\circ \pi)(T^ix)$$ and the latter is uniformly bounded by our assumption on $f \circ \pi$ and \cite[Theorem 4.11]{Gottschalk Hedlund}.
\end{proof}

\begin{lemma}\label{coboundary over any group rotation}

 Let $(G,m_{G},T_\alpha)$ be a minimal rotation on an infinite compact abelian group such that $\widehat{G}$ contains a non-torsion element. Then there exists $h\colon G \rightarrow \mathbb{R}$ such that $rh$ is 
 an $\mathbb{R}$-valued measurable but not continuous coboundary over $(G,m_{G},T_\alpha)$ for all nonzero $r \in \mathbb{R}$. 
    
\end{lemma}

\begin{proof}

    Since $\widehat{G}$ is infinite and contains a non-torsion element, it contains an infinite cyclic subgroup $C:=\{\gamma^{n}\colon n\in \mathbb{Z}\}$ for some nontrivial $\gamma \in \widehat{G}$. Then $\gamma\colon G \rightarrow \mathbb{T}$ is a topological factor map between the systems $(G,m_{G},T_{\alpha})$ and $(\mathbb{T},m_{\mathbb{T}},T_{\gamma(\alpha)})$. 
    Notice $\gamma(\alpha)$ must have infinite order in $\mathbb{T}$ for otherwise if $\gamma(n\alpha)=1$ for some $n\in \mathbb{Z}$, then $\gamma^{n}(\alpha)=1$, implying $C$ is finite. Choose $\{n_k\}$ distinct with $d(n_k\gamma(\alpha), 0)<\frac{1}{k}$ and let $n_{-k}=-n_k$ for $k\neq 0$. Defining $h\colon \mathbb{T} \rightarrow \mathbb{R}$ by $$h(x)=\sum_{k\neq 0}\frac{e(n_k\gamma(\alpha))-1}{|k|}e(n_k x),$$ we see $rh$ is a measurable but not continuous coboundary over $(\mathbb{T},m_{\mathbb{T}},T_{\gamma(\alpha)})$ for all $r \neq 0$ (cf. \cite{Furstenberg} or \cite[Lemma 3.1]{my paper}). Lemma \ref{coboundary lifts} implies $rh\circ \gamma$ is a measurable but not continuous $\mathbb{R}$-valued coboundary over $(G,m_{G},T_{\alpha})$. \end{proof}

The following lemma is used to pass from an $\mathbb{R}$-valued to $\mathbb{R}/\mathbb{Z}$-valued measurable but not continuous coboundary over any minimal system.

\begin{lemma}\label{getting a torus valued coboundary}
    Let $(X,\mu,T)$ be minimal, $q\colon\mathbb{R} \rightarrow \mathbb{R}/\mathbb{Z}$ denote the quotient map, and $nf\colon X \rightarrow \mathbb{R}$ be a measurable but not continuous coboundary over $(X,\mu,T)$ for all $n\in \mathbb{N}$. There exists $t \in \mathbb{R}$, $0< t\leq 1/4$, such that $q(tnf)$ is a measurable but not continuous coboundary for all $n\in \mathbb{N}$. 
\end{lemma}

\begin{proof}
    We use that $C(X,\mathbb{Z})$ and $C(X,\mathbb{T})/q (C(X,\mathbb{R}))$ are countable. Note that countability of the latter follows from $q (C(X,\mathbb{R}))$ being an open subgroup of the separable space $C(X,\mathbb{T})$ equipped with the uniform metric. (Suppose $f \in q(C(X,\mathbb{R}))$, let $F$ denote its continuous lift, and let $g$ be such that $\displaystyle \sup_{x\in X}d(f(x),g(x))<1/4$. Then $f(x)-g(x) \in [-1/4,1/4]$ thus it has a continuous lift which we denote by $H$. It follows $F-H$ is a continuous lift of $g$.)
    Let $$E_f=\{t \in \mathbb{R}\colon q(tf)=F_t\circ T-F_t \mathrm{\;for \;some\;continuous \;}F_t\in C(X,\mathbb{T})\}.$$ We define $p\colon  E_f \rightarrow C(X,\mathbb{T})/q(C(X,\mathbb{R}))$ by $p(t)=[F_t].$ Note $p$ is well-defined. (If $F_t'\in C(X,\mathbb{T})$ is such that $q(tf)=F_t'\circ T-F_t '$, then $F_t-F_t'$ is $T$-invariant, hence constant by minimality of $T$. Thus, $F_t-F_t'\in q(C(X,\mathbb{R}))$.) Since $C(X,\mathbb{T})/q(C(X,\mathbb{R}))$ is countable, there are countably many equivalence classes of solutions $[F_t]$. It remains to check each equivalence class is countable. 
    
    Fix $t'\in E_f$ and consider all $t$ such that $p(t')=p(t)$. It follows $F_{t'}-F_t\in q(C(X,\mathbb{R}))$. Let $G_t\in C(X,\mathbb{R})$ be such that $q(G_t)=F_{t'}-F_t$. On the other hand, $$q((t'-t)f)=q(t'f)-q(tf)=(F_t'-F_t)\circ T-(F_{t}'-F_t)=q(G_t\circ T-G_t).$$ Then $(t'-t)f-(G_t \circ T-G_t)\in C(X,\mathbb{Z})$. 
    Furthermore, we have that the map $\phi\colon p^{-1}(p(t'))\rightarrow C(X,\mathbb{Z})$ given by $$t\mapsto (t'-t)f-(G_t\circ T -G_t)$$ injects into a countable set, hence $p^{-1}(p(t'))$ is countable, therefore so is $E_f$. (If $\phi(t_1)=\phi(t_2)$, then $(t_2-t_1)f=(G_{t_1}-G_{t_2}) \circ T- (G_{t_1}-G_{t_2})$ and if $t_1\neq t_2$ then we see $$f=\frac{1}{t_2-t_1}((G_{t_1}-G_{t_2}) \circ T- (G_{t_1}-G_{t_2}))$$ is a continuous coboundary, contradicting our assumption.) We choose $t\in (0,1/4) \cap (\bigcup_{n\in \mathbb{N}}{E_{nf}})^{c}$ to conclude the proof. \end{proof}

\begin{lemma}\label{coboundaries over any nilsystem}
    Let $(X,\mu,T)$ be a minimal system with an infinite maximal equicontinuous factor $X_1$ and $\pi\colon X \rightarrow X_1$ denote a topological factor map. If $nh$ is a measurable but not continuous coboundary over $X_1$ for all $n \in \mathbb{N}$ where $h\colon X \rightarrow \mathbb{R}$, then there exists $0<t\leq 1/4$ such that for all $n\in \mathbb{N}$ we have $q(tn(h\circ\pi))$ is a measurable but not continuous $\mathbb{T}$-valued coboundary over $(X,\mu,T)$ where $q\colon \mathbb{R} \rightarrow \mathbb{R}/\mathbb{Z}$ denotes the standard projection map.
\end{lemma}

\begin{proof}
   Let $\pi\colon X \rightarrow X_1$ be a topological factor map and $nh$ be a measurable but not continuous coboundary over $X_1$ for all $n\in \mathbb{N}$. (Note such coboundaries exist by Lemmas \ref{coboundary over group rotation} and \ref{coboundary over any group rotation}.) 
    By Lemma \ref{coboundary lifts}, $nh \circ \pi$ is a measurable but not continuous coboundary over $(X,\mu,T)$. Lemma \ref{getting a torus valued coboundary} concludes the proof. 
\end{proof}

\begin{lemma} \label{torsion group sequence}
    Let $G$ be an infinite compact abelian group where $\widehat{G}$ is a torsion group and $\alpha \in G$ be such that $(G,T_{\alpha})$ is minimal. There exists a sequence $\{n_k\}$ and $\gamma_k\in \widehat{G}$ with $\mathrm{ord}(\gamma_k)=n_k$ such that $||\gamma_{_k}(\alpha)|| =1/n_k$, $d(n_k\alpha,e)\rightarrow 0$, $n_k|n_{k+1}$, and $n_k\geq 2n_{k-1}$  for all $k\in \mathbb{N}$.
\end{lemma}

\begin{proof}
    Since $G$ is monothetic, we know $\widehat{G} \leq \mathbb{S}^{1}$. Let $\alpha$ be a dense generator of $G$. First, we classify all torsion subgroups of $\mathbb{S}^1$. Let $M \leq \mathbb{S}^{1}$ be a torsion group. 
    By the Chinese remainder theorem, we can decompose $\displaystyle M= \bigoplus_{p\;prime} M_{p}$ where $M_p:=\{e\}\cup \{x\in \mathbb{S}^{1}\colon \mathrm{ord}(x)=p^n, n \geq 1\}\cap M.$ 
    Notice $M_p\leq C_{p^{\infty}}$ and all nontrivial subgroups of $C_{p^\infty}$ are $C_{p^n}$ where $1 \leq n \leq \infty$ (see, for example, \cite[Theorem 3.1]{Kaplansky}); therefore, $\displaystyle M=\bigoplus_{k=1}^{N} C_{p_k^{a_k}}$ where $N<\infty$ or $N=\infty$ for some strictly increasing sequence of primes $(p_k)$ and a sequence $(a_k)$ with $1\leq a_k \leq \infty$. 
    Using that the dual of a direct sum is a direct product we have $\displaystyle \widehat{M}=\prod_{k=1}^{N} \widehat{C_{p_k^{a_k}}}$. For each $i\in \mathbb{N}$, define $s_i^{(k)}:=\min\{k,a_i\}$. In the case $N<\infty$, define $n_k=p_1^{s_1^{(k)}}p_2^{{s_2}^{(k)}}\cdots p_N^{s_N ^{(k)}}$, so $\displaystyle\mathbb{Z}/(n_k) \mathbb{Z}\cong \bigoplus_{i=1}^{N} \mathbb{Z}/p_i^{s_i^{(k)}}\mathbb{Z}$. In the case $N=\infty$, define $n_k=p_1^{s_1^{(k)}}p_2^{{s_2}^{(k)}}\cdots p_k^{s_k ^{(k)}}$, so $\displaystyle\mathbb{Z}/(n_k) \mathbb{Z}\cong \bigoplus_{i=1}^{k} \mathbb{Z}/p_i^{s_i^{(k)}}\mathbb{Z}$. In both cases, it is clear $n_k|n_{k+1}$, $n_{k}\geq 2n_{k-1}$, and $d(n_k \alpha,e)\rightarrow 0$. Choose a sequence of characters $\phi_{k}\colon  \widehat{M} \rightarrow \mathbb{S}^{1}$ as $\phi_k(x)=e(\pi_{k}(x)/n_k)$ where $\pi_k\colon \widehat{M}\rightarrow \mathbb{Z}/(n_k) \mathbb{Z}$ is the canonical projection map. By construction $\mathrm{ord}(\phi_k)=n_k$ and for each $k$ there exists $0\leq r \leq n_k$ such that $||\phi_{k}^{r}(\alpha)||=1/n_k$. Hence, we choose $\gamma_{k}:=\phi_{k}^{r}$.
\end{proof}

\begin{lemma}\label{rotation}
    Let $G$ be a compact abelian group, $\gamma \in \widehat{G}$, $\alpha \in G$ be such that $(G,T_{\alpha})$ is minimal. There exist sequences $N_k,m_k \nearrow \infty$ such that for all $k\in \mathbb{N}$, $d(N_k \alpha, e) <1/k$, $||m_k \gamma(\alpha)||<\frac{1}{2^kN_k}$ and for all $i \in \mathbb{N}$ we have $||m_k N_i \gamma(\alpha)||\leq \frac{1}{2^k}$ for all $k \leq i$. Furthermore, defining $m_{-k}=-m_{k}$, for each $i$ we have $|f_{N_i}(x)|<1/6$, where $f\colon \mathbb{T} \rightarrow \mathbb{R}$ is given by $$f(x)=\frac{1}{48\pi}\sum_{k\neq 0} \frac{e(m_k \gamma(\alpha))-1}{|k|}e(m_k x).$$ 
\end{lemma}

\begin{proof}

 We inductively construct the sequences $m_k$ and $N_k$. Suppose for all $1\leq j \leq i-1$ we have $d(N_j \alpha, e)<1/j$ and $||m_j N_j \gamma(\alpha)||<2^{-j}$. Using that $$\overline{\Big\{(\alpha,m_1\gamma(\alpha),\dots, m_{i-1}\gamma(\alpha))^{n}: n\in \mathbb{N}\Big\}}$$ contains $(e,\textbf{0})\in G \times \mathbb{T}^{i-1}$, choose $N_i> N_{i-1
 }$ such that $d(N_i\alpha,e)<1/{i}$ and $||m_j N_i\gamma(\alpha)||< 1/2^j$ for all $1 \leq j \leq i-1$. 
 Then choose $m_i>m_{i-1}$ such that $||m_i \gamma(\alpha)||< \frac{1}{2^iN_i}$. It follows $||m_i N_i \gamma(\alpha)||\leq N_i ||m_i \gamma(\alpha)||<1/2^i$.
 

 Let $x \in \mathbb{T}$ and $f$ be defined as in the statement above. We have 
    \begin{align*}
        48\pi|f_{N_i}(x)| &\leq \sum_{k\neq 0}\Big|\frac{e(m_k N_i \gamma(\alpha))-1}{k}\Big|\leq 4\pi\Big(\sum_{0\leq k\leq i} \frac{||m_k N_i \gamma(\alpha)||}{k} +\sum_{k>i} \frac{||m_k N_i \gamma(\alpha)||}{k} \Big)\\ 
        &\leq 4\pi \Big(\sum_{0\leq k\leq i} \frac{1}{ 2^k k} + \sum_{k>i} \frac{N_i}{N_k2^k k}\Big) \leq 4\pi \Big(\sum_{0\leq k\leq i} \frac{1}{2^k} +\sum_{k>i} \frac{1}{2^k} \Big) \leq 8\pi.
    \end{align*}
    Thus $|f_{n_i}(x)|\leq 1/6$. \end{proof}


\begin{lemma}\label{beta rationally independent}
    Let $(X,\mu,T)$ be an ergodic system. For any sequence $\{n_k\}_{k \in \mathbb{N}} \subseteq \mathbb{N}$ there exists $\beta$ with $n\beta \notin \mathrm{Eig}(X,T)$ for all $n\in \mathbb{N}$ such that for infinitely many $k$, $||n_k\beta||>1/3$. 
\end{lemma}

\begin{proof}
    Since $\mathrm{Eig}(X,T)$ is countable, $m_{\mathbb{T}}\{\beta \in \mathbb{T}\colon  n\beta \notin \mathrm{Eig}(X,T)\}=1$. For completeness, we repeat the proof of \cite[Lemma 3.15]{Florian et al}. Let $B_k:=\{\beta \in \mathbb{T}\colon ||n_k \beta||\in (1/3,1/2]\}$. Since $x\mapsto n_kx$ preserves $m_{\mathbb{T}}$ for all $k\in \mathbb{N}$, it follows $m_{\mathbb{T}}(B_k)=1/3$. By Fubini's theorem, we obtain $$m_{\mathbb{T}}\{\beta \in \mathbb{T}\colon \beta \in B_k\;\mathrm{for \;infinitely\; many\; k\in \mathbb{N}}\}\geq 1/3.$$ The intersection of these two sets has positive measure, in particular, is non-empty. \end{proof}

\section{Recurrence not lifting over nilsystems}

We begin this section with a proposition which is key for showing property 1 in Theorem \ref{main theorem} and may be of independent interest. Its essence was partially touched upon in the proof of \cite[Proposition 4.4]{my paper} and we repeat many of the same techniques here. 

\begin{proposition}\label{the should be true lemma}
Let $(X,\mu,T)$ be a strictly ergodic system 
such that $Z_k(T)\overset{m}\cong X_k(T)$ for all $k\in \mathbb{N}$. Let also $h\colon X \rightarrow \mathbb{T}$ be continuous, $nh$ be a measurable but not continuous coboundary over $(X,\mu,T)$ and $n\beta \notin \mathrm{Eig}(X,T)$ for all $n \in \mathbb{N}$. Then $$T_{h+\beta}(x,y)=(Tx, y+h(x)+\beta), \;\; T_{h+\beta}\colon X \times \mathbb{T} \rightarrow X \times \mathbb{T},$$ is strictly ergodic, and for all $k \in \mathbb{N}$ we have $X_k(T)\overset{t}\cong X_k(T_{h+\beta})$ and $Z_1(T_{h+\beta})\overset{m}\ncong X_1(T_{h+\beta})$.
\end{proposition}

\begin{proof}
  We observe if $h=G\circ T-G$ for some measurable but not continuous $G\colon X \rightarrow \mathbb{T}$, then the map $\psi\colon X\times \mathbb{T} \rightarrow X \times \mathbb{T}$ where
    $$\psi(x,y)=(x,y-G(x))$$ is a measurable isomorphism between the systems $(X\times \mathbb{T},\mu \times m_{\mathbb{T}},T_{h+\beta})$ and $(X\times \mathbb{T}, \mu \times m_{\mathbb{T}},T \times T_{\beta})$. 

    It is clear $T_{h+\beta}$ preserves $\mu \times m_{\mathbb{T}}$. Since $n\beta \notin \mathrm{Eig}(X,T)$, $T$ and $T_{\beta}$ share no common nontrivial eigenvalues other than 1; hence $T \times T_{\beta}$ is ergodic. Using the measurable isomorphism $\psi$, we find that $T_{h+\beta}$ must also be ergodic. To obtain strict ergodicity, we draw on classical results of Furstenberg \cite{Fur81}(see proof of \cite[Proposition 4.4]{my paper} for details).
    
    We recall the elementary fact that given ergodic systems $(Y,m,S)$, $(Y',m',S')$ and $(Y\times Y',m\times m',S\times S')$ we have $$Z_k(S\times S')\overset{m}\cong (Z_k(Y)\times Z_k(Y'),m_{Z_k(Y)}\times m'_{Z_k(Y')},S\times S').$$ If we additionally use our measurable isomorphism $\psi$ from earlier and that $Z_k(T)\overset{m}\cong X_k(T)$, we obtain 
    \begin{align*}
    Z_k(T_{h+\beta})&\overset{m}\cong Z_k(T\times T_{\beta})\overset{m}\cong  (Z_k(X)\times \mathbb{T},m_{Z_k(X)}\times m_{\mathbb{T}},T \times T_{\beta})\overset{m}{\cong}\\&(X_k(X)\times \mathbb{T},m_{X_k(X)}\times m_{\mathbb{T}},T \times T_{\beta}).
    \end{align*}
Let $\phi\colon Z_k(T_{h+\beta})\rightarrow (X_k(X) \times \mathbb{T},m_{X_k(X)}\times m_{\mathbb{T}},T \times T_{\beta})$ denote this measurable isomorphism. Since $(X_k(X) \times \mathbb{T},m_{X_k(X)}\times m_{\mathbb{T}},T \times T_{\beta})$ has $X_k(T)$ as a factor, we have a measurable factor map $\tau_1\colon Z_k(T_{h+\beta}) \rightarrow X_k(T)$. Furthermore, we have a measurable factor map $\tau_2\colon Z_k(T_{h+\beta}) \rightarrow X_k(T_{h+\beta})$. It follows by the construction of the measurable nilfactors $Z_k$ in \cite[Chapter 13, Proposition 16]{HK2} that there exists a measurable factor map $\tau_3\colon  X_k(T_{h+\beta}) \rightarrow X_k(T)$ such that $\tau_1=\tau_3\circ \tau_2$. We illustrate this in the commutative diagram below.
\[
\begin{tikzcd}[column sep=small, row sep=small]
& Z_k(T_{h+\beta}) \arrow[dl, "\tau_1"'] \arrow[dr, "\tau_2"] & \\
X_k(T) & & X_k(T_{h+\beta}) \arrow[ll, "\tau_3"', dashed]
\end{tikzcd}
\]
Then $\tau_1\circ \phi^{-1}=\tau_3\circ \tau_2\circ \phi^{-1}$.  Thus, $X_k(T_{h+\beta})$ lies as an intermediate measurable factor between $(X_k(X) \times \mathbb{T},m_{X_k(X)}\times m_{\mathbb{T}},T \times T_{\beta})$ and $X_k(T)$.
   Proceeding via the same analysis as in Case 1 of \cite[Proposition 4.4]{my paper}, the nontrivial intermediate factors that can give rise to $X_k(T_{h+\beta})$, up to measurable isomorphism, are $$(X_k(X)\times \mathbb{T}, m_{X_k(X)} \times m_{\mathbb{T}},T\times T_{n\beta}),\;n\in \mathbb{N}.$$ Since the latter is also an inverse limit of nilsystems, by \cite[Theorem A.1]{topological nilfactors}, this measurable isomorphism agrees almost everywhere with a topological isomorphism. By \cite[Lemma 4.1]{my paper}, it suffices to show that for all $n \in \mathbb{N}$ there are no topological factor maps to $(X_k(X)\times \mathbb{T},m_{X_k(X)}\times m_{\mathbb{T}},T\times T_{n\beta})$ to obtain $X_k(T)\overset{t}\cong X_k(T_{h+\beta})$.

    Any such continuous map $\pi\colon X \times \mathbb{T} \rightarrow X_k(X) \times \mathbb{T}$ is of the form 
    $$\pi(x,y)=(\pi_1(x,y),\pi_2(x,y))$$ for some $\pi_1: X \times \mathbb{T} \rightarrow X_k(X)$ and $\pi_2: X \times \mathbb{T} \rightarrow \mathbb{T} $ continuous. Suppose $\pi$ was additionally a topological factor map. Then  $$(\pi_1\circ T_{h+\beta},\pi_2 \circ T_{h+\beta})=(T\circ \pi_1, T_{n\beta}\circ \pi_2 ).$$ By equating the second coordinates we obtain a coboundary equation 
    \begin{equation}\label{equating second coordinates}
        \pi_2 \circ T_{h+\beta}-\pi_2=n\beta.
    \end{equation}
    Notice $\pi_2(x,y)=ny-nG(x)$ is a solution to (\ref{equating second coordinates}). Since $T_{h+\beta}$ is ergodic, all solutions $\pi_2$ must be unique up to a constant almost everywhere, hence of the form $ny-nG(x)+c$ for some $c \in \mathbb{T}$. By Fubini's theorem, there exists $y_0\in \mathbb{T}$ such that for $\mu$-a.e. $x$, 
    \begin{equation}\label{fubini}
    nG(x)=-\pi_2(x,y_0)+ny_0+c.
    \end{equation}
    Since (\ref{fubini}) holds $\mu$-almost everywhere, $\supp(\mu)=X$, and the right-hand side of (\ref{fubini}) is continuous, we obtain (\ref{fubini}) holding for every $x$, hence $nG$ is continuous. This contradicts our assumption that $nh=(nG)\circ T-nG$ is not a continuous coboundary.

    To verify $Z_1(T_{h+\beta}) \overset{m}\ncong X_1(T_{h+\beta})$, it suffices to find an eigenfunction of $T_{h+\beta}$ which is not continuous. Indeed, $\eta(x,y)=y-G(x), \eta\colon X \times \mathbb{T}\rightarrow \mathbb{T}$, is an eigenfunction with eigenvalue $\beta$ which is not continuous since $G$ is not continuous. \end{proof}

    In the topological proof of \cite[Theorem 4.1]{Host Kra Maass} in \cite[Section 4.3]{Host Kra Maass}, it is shown every set of recurrence for all finite dimensional rotations is a set of recurrence for any nilsystem. The proof provides the following stronger statement: If $(X,T)$ is a minimal nilsystem and if $R \subseteq \mathbb{N}$ is a set of recurrence for $(X_1,T)$, then it is a set of recurrence for $(X,T)$. This is directly seen in statement (4) of their topological proof.




\subsection{\textbf{Proof of Theorem \ref{main theorem}}}\label{proof of main theorem}

Let $G=\overline{\langle \alpha \rangle}$ and $(G,m_{G},T_{\alpha})$ denote the $X_1$ factor of a minimal nilsystem $(X,\mu,T)$ where $\pi\colon X \rightarrow G $ is a topological factor map. Let also $q\colon \mathbb{R} \rightarrow \mathbb{R}/\mathbb{Z}$ denote the standard quotient map. We separate the proof into the case when the dual group $\widehat{G}$ is a torsion group and when it is not.

\textit{Case 1. $\widehat{G}$ is a torsion group}

By Lemma \ref{torsion group sequence}, there exists a sequence $\{n_k\}$ and $\gamma_k\in \widehat{G}$ with $\mathrm{ord}(\gamma_k)=n_k$ such that $||\gamma_{_k}(\alpha)|| =1/n_k$, $d(n_k\alpha,e)\rightarrow 0$, $n_k|n_{k+1}$, and $n_{k}\geq 2n_{{k-1}}$ for all $k\in \mathbb{N}$. Define $n_{-k}=n_k$ for all $k\in \mathbb{N}$. It follows from Lemma \ref{beta rationally independent} that there exists $\beta \in \mathbb{T}$ with $n\beta \notin \mathrm{Eig}(X,T)$ and there exists a subsequence $\{n_{s_k}\}$ such that $||n_{s_k}\beta||>1/3$. 
Define $h\colon G \rightarrow \mathbb{R}$ by $$h(g)=\frac{1}{6\pi}\sum_{k\neq 0} \frac{\gamma_{k}(\alpha)-1}{|k|}\gamma_{k}(g)$$ where $\gamma_{-k}:=\overline{\gamma_k}$. 
Since in particular $||\gamma_{k}(\alpha)||<1/k$, the construction in Lemma \ref{coboundary over group rotation} shows $nh$ is a measurable but not continuous coboundary over $(G,m_G,T_{\alpha})$ for all $n\in \mathbb{N}$. By Lemma \ref{coboundaries over any nilsystem}, there exists $0<t\leq 1/4$ such that $q(t(nh\circ \pi))$ is a measurable but not continuous coboundary over $(X,\mu,T)$. 
Choosing $$f=q(t(h\circ \pi))+\beta$$ we see property 1 is satisfied by Proposition \ref{the should be true lemma}. 
Recall that $n_{s_k} \alpha \rightarrow e$, hence $$R=\{n_{s_i}\colon i \in \mathbb{N}\}$$ is a set of recurrence for $(G,m_G,T_{\alpha})$. Furthermore, it is a set of recurrence for $(X,\mu,T)$ by \cite[Theorem 4.1]{Host Kra Maass}. To show $R$ is not a set of recurrence for $(X \times \mathbb{T},\mu\times m_{\mathbb{T}},T_f)$, it suffices to show for all $x\in X$ and $n_{s_i} \in R$, $$ \Big|\Big|\sum_{k=0}^{n_{s_i} -1} th(\pi(T^kx))+n_{s_i}\beta\Big|\Big|\geq 1/6.$$ (Notice $\displaystyle T_{f}^n(x,y)=\Big(T^nx,y+\sum_{k=0}^{n-1} t(h\circ \pi) (T^kx)+n\beta\Big)$.) We first observe the following
\begin{align*}
    6\pi|t(h\circ \pi)_{n_{i}}(x)|&=6\pi\Big|\sum_{k=0}^{n_i -1} th(\pi(x)+k\alpha)\Big| \leq t\sum_{k\neq 0}\frac{|\gamma_{k}(n_i\alpha)-1|}{k}\\&\leq 4t\pi\Big(\sum_{0\leq k\leq i}\frac{||\gamma_{k}(n_i\alpha)||}{k} +\sum_{k>i}\frac{||\gamma_{k}(n_i \alpha)||}{k}\Big) 
\end{align*}
Since $n_k|n_{k+1}$ for all $k\in \mathbb{N}$, we have $\gamma_{k}(n_i\alpha)=\gamma_{k}^{n_i} (\alpha)=1$ whenever $0\leq k \leq i$. Therefore, $\displaystyle \sum_{0\leq k\leq i}\frac{||\gamma_{k}(n_i\alpha)||}{k} =0$. 
On the other hand, $$4t\pi\sum_{k>i} \frac{||\gamma_{k}(n_i \alpha)||}{k} \leq 4t\pi\sum_{k>i}\frac{n_i ||\gamma_{k}(\alpha)||}{k}=4t\pi \sum_{k> i} \frac{n_i}{n_kk}\leq 4t\pi\sum_{k> i}\frac{n_i}{n_k} \leq \pi$$ therefore $|t(h\circ \pi)_{n_i}|$, and in particular $|t(h\circ \pi)_{n_{s_i}}|$, are uniformly bounded by 1/6. Finally, because $||n_{s_i}\beta||>1/3$ for each $i$ we obtain $$\Big|\Big|\sum_{k=0}^{n_{s_i} -1} th(\pi(T^kx))+n_{s_i}\beta\Big|\Big| \geq 1/3-1/6=1/6.$$

\textit{Case 2. $\widehat{G}$ contains an infinite cyclic subgroup}

Let $\{\gamma^n\colon n\in \mathbb{Z}\}\leq \widehat{G}$ be our infinite cyclic subgroup. Employing Lemma \ref{rotation}, we construct sequences $m_k, N_k \nearrow \infty$ with $||m_k \gamma(\alpha)||<\frac{1}{2^k N_k}$, $d(N_k\alpha,e)<1/k$, $||m_k N_i \gamma(\alpha)||<\frac{1}{2^k}$ for all $k\leq i$ and $m_{-k}=-m_k$. By passing to a subsequence $\{N_{s_k}\}$, Lemma \ref{beta rationally independent} provides $\beta\in \mathbb{T}$ with $n\beta \notin \mathrm{Eig}(X,T)$ for all $n\in\mathbb{N}$ such that $||N_{s_k}\beta||>1/3$. 
By construction, $$R=\{N_{s_i}\colon i\in \mathbb{N}\}$$ is a set of recurrence for $(G,m_G,T_{\alpha})$. Furthermore, $R$ is also a set of recurrence for $(X,\mu,T)$ by \cite[Theorem 4.1]{Host Kra Maass}. Since $(X,\mu,T)$ is minimal, its factor $(G,m_{G},T_{\alpha})$ must also be minimal. 
Then Lemma \ref{coboundary over any group rotation} shows $nh\circ \gamma$ is a measurable but not continuous coboundary over $(G,m_{G},T_{\alpha})$ for all $n\in \mathbb{N}$ where recall $$h(x)=\frac{1}{48\pi}\sum_{k\neq 0} \frac{e(m_k \gamma(\alpha))-1}{|k|}e(m_k x).$$ By Lemma \ref{coboundaries over any nilsystem}, there exists $0<t\leq 1/4$ such that $q(t(nh \circ \gamma \circ \pi))$ is a measurable but not continuous coboundary over $(X,\mu,T)$ for all $n\in \mathbb{N}$. Choosing $f=q(th\circ \gamma \circ \pi)+\beta$ concludes property 1 using Proposition \ref{the should be true lemma}. 
By Lemma \ref{rotation}, for each $i\in \mathbb{N}$ and $x \in X$, $|h_{N_{s_i}}(x)|$ is uniformly bounded by 1/6, hence $|(th\circ \gamma\circ \pi)_{N_{s_i}}(x)|<t/6$. Since $||N_{s_i}\beta||>1/3$, we observe $$\displaystyle \Big|\Big|\sum_{k=0}^{N_{s_i} -1} (th\circ \gamma\circ \pi)(T^kx))+N_{s_i}\beta\Big|\Big|\geq 1/3-t/6\geq7/24,$$ hence $R$ is not a set of recurrence for $(X \times \mathbb{T},\mu \times m_{\mathbb{T}},T_{f})$.\;\;\;\;\;\;\;\;\;\;\;\;\;\;\;\;\;\;\;\;\;\;\;\;\;\;\;\;\;\;\;\;\;\;\;\;\;\;\;\;\;\;\;\;\;\;\;\;\;\;\;\;\;\;\;\;\;\;\;\;\;\;\;\;\;\;\;\;\;\;\;\;\;\;\;\;\;\;\;\;\;\;\;\;\;\;\;\;\;\;\;\;\;\;\;\;\;\;\;$\square$

\section{Hidden frequencies for some non-distal systems}

The systems we explored in Theorem \ref{main theorem} were all distal, meaning for any $x\neq y$ in $X$, $\displaystyle\inf_{n\in \mathbb{N}}d(T^n x, T^n y)>0$. Rotations are the simplest examples of distal systems. Systems exhibiting some random behavior in their dynamics are topologically weak mixing, 
that is, its maximal equicontinuous factor is trivial. It is tempting to expect that a set of recurrence for $(\mathbb{T},T_{\alpha})$ would be a set of recurrence for $(\mathbb{T} \times X, T_{\alpha} \times T)$ where $(X,T)$ is topologically weak mixing. However, we show this does not always hold. 

\begin{proposition}\label{top weak mixing}
    For every minimal topologically weak mixing but not topologically mixing system $(X,T)$ there exists a set $R \subseteq \mathbb{N}$ and $\alpha \in \mathbb{T}\setminus\mathbb{Q}$ such that $R$ is a set of recurrence for $(\mathbb{T},T_{\alpha})$ but it is not a set of recurrence for $(\mathbb{T}\times X,T_{\alpha}\times T).$
\end{proposition}

The next lemma follows immediately from definitions. Note a system $(X,T)$ is \textit{topologically mixing} if for any two open sets $U,V \subseteq X$ there exists $N\in \mathbb{N}$ such that for all $n \geq N$, we have $U\cap T^{-n} V \neq \emptyset$.

\begin{lemma}\label{transitive but not mizing}
    Let $(X,T)$ be transitive but not topologically mixing. Then there exists an open set $U \subseteq X$ and a sequence $n_k \nearrow \infty$ such that $U \cap T^{-n_k}U=\emptyset$ for all $k \in \mathbb{N}$. 
\end{lemma}

\begin{proof}
   We show that if for every nonempty open set $W \subseteq X$ there exists $N\in \mathbb{N}$ such that for all $n\geq N$, $W \cap T^{-n}W\neq \emptyset$, then for all open $U,V \subseteq X$ there exists $M\in \mathbb{N}$ such that for all $n\geq M$, $U \cap T^{-n}V\neq \emptyset.$ By transitivity, there exists $n\in \mathbb{N}$ such that $U \cap T^{-n} V \neq \emptyset$. Since the latter is open, there exists $M \in \mathbb{N}$ such that for all $m \geq M$ we have $U \cap T^{-n} V \cap T^{-m}U \cap T^{-(m+n)}V \neq \emptyset$. Then for all $j \geq M+n$ we have $U \cap T^{-j}V\neq \emptyset.$ \end{proof}


\begin{lemma}\label{set of recurrence not lifting for wm}
    For any infinite set $R\subseteq \mathbb{N}$ there exists $\alpha \in \mathbb{T} \setminus{\mathbb{Q}}$ such that $R$ is a set of recurrence for $(\mathbb{T},T_{\alpha}).$
\end{lemma}

\begin{proof}
     Let $\displaystyle I_{k}:=\bigcup_{r \in R}\{\alpha \in \mathbb{T}\colon ||r\alpha||<1/k\}.$ For any open interval $I$, there exists $j\in \mathbb{N}$ and $r \in R$ such that $j/r \in I$. Letting $\beta=j/r$ we see $\beta \in I_k$, so $I_k$ are open dense sets for each $k$. By the Baire category theorem, $\displaystyle\bigcap_{k\in \mathbb{N}} I_{k}$ is comeager, hence contains an irrational.
\end{proof}

Let $(X,T)$ be any minimal but not topologically mixing system and $\{n_k\}$ be the sequence obtained from Lemma \ref{transitive but not mizing}. By choosing $R=\{n_k\colon k \in \mathbb{N}\}$ in Lemma \ref{set of recurrence not lifting for wm}, we obtain Proposition \ref{top weak mixing}. Here, we use that in a minimal system if $R\subseteq \mathbb{N}$ is a set of recurrence for $(X,T)$, then for every open set $U\subseteq X$ there exists $n\in R$ such that $U \cap T^{-n}U \neq \emptyset$. Indeed, let $U\subseteq X$ and $x\in X$ be such that $\displaystyle\inf_{n\in R}d(x,T^n x)=0$. By minimality of $(X,T)$, there exists $j\in \mathbb{N}$ such that $T^j x \in U$. By continuity of $T$, $\displaystyle\inf_{n\in R}d(T^jx,T^{n+j} x)=0$. Hence there exists $n \in R$ such that $T^jx\in U \cap T^{-n}U$.

\end{document}